\documentclass[10pt]{amsart}

\usepackage{amsfonts}

\usepackage{graphicx}
\usepackage{amsmath}
\usepackage{setspace}
\usepackage[tmargin=2.0cm, bmargin=2.0cm, lmargin=3.0cm, rmargin=3.0cm]{geometry}

\def\BEN{\begin{enumerate}}  \def\BI{\begin{itemize}}
\def\EEN{\end{enumerate}}   \def\EI{\end{itemize}}
   \def\sec{\section} 

\def\beq{\begin{eqnarray}} \def\eeq{\end{eqnarray}}

\def\al*#1{\begin{align*}#1\end{align*}}

\def\ga*#1{\begin{gather*}#1\end{gather*}}

\def\alat*#1#2{\begin{alignat*}{#1}#2\end{alignat*}}
\def\bea{\begin{eqnarray*}}
\def\eea{\end{eqnarray*}}
\def\ml*#1{\begin{multline*}#1\end{multline*}}

\def\P{{\mathbb P}}   
\def\E{{\mathbb E}}   
\def\N{{\mathbb N}}

\newtheorem{Thm}{Theorem}
\newtheorem{Fa}{Lemma}

 \newtheorem{Def}{Definition}
\newtheorem{Rem}{Remark} 
 
\newcommand{\exit}{{\mbox{\, \vspace{3mm}}}
\hfill\mbox{$\square$}}

\title{Discrete time ruin probability \\ with Parisian delay}

\author{Irmina Czarna}
\address{Mathematical Insititute, University of Wroc\l aw, pl. Grunwaldzki 2/4, 50-384 Wroc\l aw, Poland}
\email{irmina.czarna@gmail.com}

\author{Zbigniew Palmowski}
\address{Department of Applied Mathematics, Wroc\l aw University of Science and Technology, ul. Wyb. Wyspia\'nskiego 27, 50-370 Wroc\l aw, Poland}
\email{zbigniew.palmowski@gmail.com}

\author{Przemys{\l}aw \'Swi\c{a}tek}
\address{Credit Suisse}
\email{przemek.swiatek@gmail.com}

\thanks{This work is partially supported by National Science Centre Grant No. DEC-2011/01/B/HS4/00982
(2012-2013). All authors kindly acknowledges partial support by the project RARE -318984, a Marie Curie IRSES Fellowship within the 7th European
Community Framework Programme.}

\date{\today}

\begin{document}

\begin{abstract}
In this paper we evaluate the probability of the discrete time Parisian ruin that occurs when surplus process stays below or at zero at least for some fixed duration of time $d>0$. We identify expressions for the ruin probabilities within finite and infinite-time horizon. We also find their light and heavy-tailed asymptotics
when initial reserves approach infinity. Finally, we calculate these probabilities for a few explicit examples.

\vspace{3mm}


\end{abstract}

\maketitle

\pagestyle{myheadings} \markboth{\sc I.\ Czarna --- Z.\ Palmowski
--- P.\ \'Swi\c{a}tek} {\sc Discrete time ruin probability with Parisian delay}

\vspace{1.8cm}

\tableofcontents

\newpage

\sec{Introduction}\label{sec:intro}
In the present paper we consider the following process:
\begin{eqnarray}\label{CL}
 R_n = u +  n -S_n,
\end{eqnarray}
where $u>0$ denotes the initial reserve and
$$S_n=\sum\limits_{i=1}^n Y_i,\qquad n=1,2,3,\dots .$$
We assume that  $Y_i$ ($i=1,2,\dots$) are i.i.d. claims
and we also assume that premium rate equals to $1$.
We denote $\P(Y_1=k)=p_k$ for $k=0, 1, 2, \ldots$ and we assume that $\mu=\E(Y_1)<1$, hence $R_n\to +\infty$ a.s.
The risk process $R$ starts from $R_0=u$ and
later we use convention $\P(\cdot|R_0=u)=\P_u(\cdot)$ and $\P_0=\P$.
The discrete-time model (\ref{CL}) is very important for actuarial practice,
since many crucial quantities related to this model have a recursive nature and are readily programmable in practice; see
e.g., \cite{Lietal, k} and references therein.

One of the most important characteristics in risk theory is finite-time ruin probability defined by
$\P_u(\tau^0<t)$  for the ruin moment
$\tau^0=\inf\{n \in \N: R_n \leq 0 \}$
and fixed time horizon $t$.
Let us note here that our definition is compatible with many papers, see e.g., Gerber \cite{h} and Dickson \cite{f}.
Other authors define the ruin moment when the reserve takes strictly negative value (see e.g., Willmot \cite{k}).
In this paper we extend this notion to so-called Parisian ruin probability,
which occurs if the process $R$ stays below or at zero at least for a fixed time period $d \in\{1,2,\ldots\}$.
Formally, we define Parisian ruin time by:
$$\tau^{d}=\inf\{n\in \N: n-\sup\{s<n: R_s > 0\}> d, R_n\leq 0 \}$$
and we consider Parisian ruin probabilities $\P_u(\tau^d<t)$ and $\P_u(\tau^d<\infty)$.

The case $d=0$ corresponds to the classical ruin problem.
There are already a number of relevant results analyzing this case, e.g.,
Dickson and Hipp \cite{13}, Gerber and Shiu \cite{ 15,16}, Li and Garrido \cite{22,23}, Lin and Willmot \cite{ 24, 25}, Shiu \cite{ j}, Willmot \cite{33}. Moreover,
Li et al. in \cite{Lietal} presented a comprehensive review.
The discrete model was first proposed by seminal paper of Gerber \cite{h}. In this paper
the ruin probability was expressed in terms of total amount $S_t$ of claims cumulated up to time $t$.
Explicit formulas for the ruin probability were also derived by Willmot \cite{k} (see also Cheng et al. \cite{a}), where the author used analytical
techniques, such as Lagrange's expansions of moment generating functions.

Other related results concern the expected discounted penalty function (so-called Gerber-Shiu function),
corresponding to the joint distribution of the surplus immediately
before and at ruin moment. For example
Cheng et al. in \cite{a} considered the discounted probability of ruin: $\sum_{n=1}^\infty \upsilon^n \P_u(R_{\tau^0-1}=x, R_{\tau^0}=-y, \tau^0=n)$, where the surplus just before ruin is $x$,
the deficit at ruin equals $y$ and $\upsilon$ is a discount factor $(0<\upsilon<1)$. Li and Garrido \cite{21} explored this topic further giving a recursive formula for the expected discounted penalty function due to ruin. In their proof they used the moment generating functions. In continuous-time
model a similar approach was applied by Dickson \cite{10}.
A detailed discussion was given when the claim size is geometrically distributed.
Another approach is based on a defective renewal equation; see Landriault \cite{17},
Pavlova and Willmot \cite{29}.
Results for the discrete-time risk models were also used as approximations or bounds for the corresponding results in continuous time,
see Cossette et al. \cite{c} and Dickson et al. \cite{12} for the approximating procedures.
Other references on the related topics are: Cossette et al. \cite{b, d}, Dickson \cite{ f, g}, Li \cite{i, 20}, Michel \cite{ 26}, Wu and Li \cite{35}, Yang et al. \cite{ 36}, Yuen and Guo \cite{ 37, 38}.

The name for the problem considered in this paper is borrowed from Parisian option, where prices
are activated or canceled depending on a type of option when
underlying asset stays above or below barrier long enough (see Albrecher et al. \cite{Hansjoerg}, Chesney et al. \cite{ chesneyetal}, Dassios and Wu \cite{DassWu1}).
We believe that Parisian ruin probability could be a better measure of risk in many situations
giving possibility for insurance company to get solvency. So far the Parisian ruin probability has been considered only
in a continuous-time setting. In particular, Dassios and Wu \cite{DassWu1} analyze the continuous-time classical risk process (\ref{CL}) with exponential claims
and the Brownian motion with drift. Dassios and Wu \cite{DassWu2} found also Cram\'{e}r-type asymptotics for this risk process.
Czarna and Palmowski \cite{IZ} and Loeffen et al. \cite{RIZ}
extended these results to the case of a general spectrally negative
L\'{e}vy process using the fluctuation theory. Another
possible way of defining the Parisian delay is based on exchanging  the deterministic, fixed delay $d$
by an independent exponential random variable; see e.g., Landriault et al. \cite{Land2} and Baurdoux et al. \cite{erik}.

The main goal of this paper is to derive discrete-time counterparts of the results above and propose efficient numerical procedure for finding Parisian probability of (non-)ruin within finite time.

This paper is organized as follows.
In Sections \ref{sec:main} and \ref{sec:infinite} we give main representation of Parisian non-ruin and ruin probabilities within finite and infinite time respectively.  In Sections \ref{sec:cram} and \ref{sec:heavy}
we give asymptotics of Parisian ruin probability in Cram\'{e}r (light-tailed) and heavy-tailed cases.
Finally, in Section \ref{sec:examples} we analyze
a few particular examples.


\sec{Parisian non-ruin probability over any finite-time horizon }\label{sec:main}

In the main result we will use the following Seal-type formula proved in Lef\'evre and Loisel \cite[Prop. 2.4]{le}.
\begin{Fa}\label{classic_ruin}
We have $\P_u(\tau^0=1)=\P(Y_1\geq u+1)$ and for $t\geq 1$:
\begin{eqnarray}
\P_u(\tau^0\geq t+1)=\sum_{j=0}^{u+t-1}p_j^{*t}-\sum_{j=u+1}^{u+t-1}p_j^{*(j-u)}\left(\sum_{k=j}^{u+t-1}\frac{t+u-k}{t+u-j}p_{k-j}^{*(t+u-j)}\right),
\end{eqnarray}
where $\{p_k^{* t}, n\in \N  \}$ denotes the $t$-th convolution of the law of $Y_1$.
\end{Fa}

\begin{Rem}
There exist alternative, recursive ways of calculating $\P_u(\tau^0\geq t+1)$. For example,
De Vylder and Goovaerts \cite{VG} give the following procedure:
\begin{equation}\label{recfr}
\P_u(\tau^0\geq 0)=1,\qquad \P_u(\tau^0\geq t)=\sum_{k=0}^u p_k\P_{u+1-k}(\tau^0\geq t-1);
\end{equation}
see also Dickson and Waters \cite{DicksonWaters}.
\end{Rem}

\begin{Fa}\label{Loisel}
For $s\geq 1$ we have:
\begin{eqnarray}\nonumber
\lefteqn{\P_{u}(\tau^{0}=s,-R_{\tau^0}=z)=\sum_{k=0}^{u+s-2}\P_{u}(\tau^0>s-1,S_{s-1}=k)p_{u+s-k+z}}\\\nonumber&=&
\sum_{k=0}^{u+s-2}p_k^{*(s-1)}p_{u+s-k+z}-
\sum_{k=u+1}^{u+s-2}\sum_{j=u+1}^{k}\frac{s-1+u-k}{s-1+u-j}p_{k-j}^{*(s-1+u-j)}p_j^{*(j-u)}p_{u+s-k+z}.\\
\end{eqnarray}
\end{Fa}
\proof
The first equality is a straightforward consequence of Markov property.
The second equality follows from
decomposition of a trajectory of $R_n$ into two parts and from lemma given in Lef\'evre and Loisel
\cite[Lem. 2.3]{le}.
\exit\\
\begin{Rem}\label{geometric} \rm
Assume that the claim $Y_1$ is of the form $Y_1=I\cdot B$, where $I$ is a Bernoulli random variable with $\P(I=1)=1-\P(I=0)=b$ and $B$ is a geometric random variable with parameter $1-q$, that is \begin{equation}\label{bingeom}
\P(Y_1=0)=p_0=1-b,\qquad \P(Y_1=z)=p_z=bq^{z-1}(1-q)\qquad {\rm for}\quad z=1,2,\ldots .
\end{equation}
Model (\ref{CL}) with such distribution of claims is a particular case of widely used compound binomial model.
In this model at each unit of time a claim might arrive with probability $b$ and there is no claim with probability $1-b$
(see Cheng et al. \cite{a}, Li and Guo \cite{ liuguo}, dos Reis \cite{egido}).
By lack of memory property of geometric distribution, we have:
\begin{equation}\label{geombern1}\nonumber
\P_{u}(\tau^{0}=s,-R_{\tau^0}=z)=\P_{u}(\tau^{0}=s)q^{z}(1-q).
\end{equation}
Similarly,
\begin{equation}\label{geombern2}\nonumber
\P_{u}(\tau^{0}<\infty,-R_{\tau^0}=z)=\P_{u}(\tau^{0}<\infty)q^{z}(1-q).
\end{equation}
Moreover, Wu and Li \cite[eq. (4.7)]{35} showed that for the compound binomial model:
\begin{equation}\label{geom1}\nonumber \P_{u}(\tau^{0}<\infty)=\xi \left[q+\xi(1-q)\right]^{(u-1)},
\end{equation}
where $\P_{1}(\tau^{0}<\infty)=\xi$ solves
$\xi=\widehat{k}\left(q+\xi(1-q)\right)$
for p.g.f $\widehat{k}(z)$ of generic claim $Y$.
In our case of (\ref{bingeom}) we have
$\widehat{k}(z)=\frac{zb}{1-z(1-b)}$
and hence
$$\xi=\frac{bq}{(1-q)(1-b)}\;.$$
\end{Rem}
Let $$\tau_x=\inf\{n \in \N: R_n=x\}.$$
The main representation of the finite-time Parisian ruin probability is given in the next theorem.
\begin{Thm}\label{ThmMain}
For $u\geq 1$, the recursive representation of the Parisian non-ruin probability until finite-time t is as follows.
For $t \leq d+1$ we have $\P_{u}(\tau^{d}\geq t)=1.$\\
For $t\geq d+2$:
\begin{eqnarray}\label{finite}\nonumber
\P_{u}(\tau^{d}\geq t)=\P_{u}(\tau^0\geq t-d) +
\sum\limits_{s=1}^{t-d-1}\sum\limits_{\omega=1}^{d}\sum\limits_{z=0}^{\omega-1}\P_{u}(\tau^0=s,-R_{\tau^0}=z)\P(\tau_{z+1}=\omega)\P_1(\tau^{d}\geq t-{\omega}-s),\\
\end{eqnarray}
where
\begin{equation}\label{Kendall}
\P(\tau_x=\omega)=\frac{x}{\omega}\P(R_{\omega}=x)=\frac{x}{\omega}p_{\omega-x}^{*\omega}
\end{equation}
 and probabilities
$\P_{u}(\tau^0\geq t)$, $\P_{u}(\tau^0=s,-R_{\tau^0}=z)$ are given in Lemmas \ref{classic_ruin} and \ref{Loisel}, respectively.
\end{Thm}
\proof
The statement for $t \leq d+1$ is obvious.
For $t\geq d+2$ the Parisian ruin occurs after time $t$ if and only if
one of the following two separate scenarios happen. In the first one, the classical ruin time $\tau^0$
happens after time $t-d$. In the second scenario,
we can decompose possible trajectory that drops below zero into two parts.
The first part runs until the first time it hits $1$. The second one
runs after this time.
Precisely, the first piece of this trajectory crosses level $0$ at ruin time $\tau^0\in\{1,\ldots, t-d-1\}$ and it has the undershoot of size $-z\leq 0$. Later it returns to level $1$. This excursion from $-z$
must have length $\omega \leq d$, otherwise we will have Parisian ruin before time $t$.
The second part of the above mentioned trajectory starts at $1$ and avoids Parisian ruin over $t-\omega-s$. This observation and the strong Markov property imply formula (\ref{finite}).
Equality (\ref{Kendall}) follows from Kendall's identity for a random walk (see Alili et al. \cite{alilietal}, Bertion \cite[Cor. VII.3]{bertbook} and Feller \cite[eq. (9.3), p. 424]{Feller}).
\exit
\begin{Rem}\rm
To compute probability given in (\ref{finite}) for $u\geq 2$ we start from counting $\P_1(\tau^{d}\geq t)$ using the algorithm given in Theorem \ref{ThmMain}
for $u=1$.
\end{Rem}

\section{Ultimate Parisian ruin probability}\label{sec:infinite}
Gerber \cite{h}, Lef\'evre and Loisel \cite[Cor. 2.8]{le}, Shiu \cite{j} and Willmot \cite{k} proved the following result.
\begin{Fa}\label{infty_Loisel}
We have:
\begin{equation}
\P_u(\tau^0<\infty)=(1-\mu)\sum_{j=u+1}^{\infty}p_{j}^{*(j-u)}.
\end{equation}
\end{Fa}

\begin{Thm}\label{ThmMain2}
For $u\geq 1$ the representation of ultimate Parisian ruin probability is given by:
\begin{eqnarray}\label{infinite}\nonumber
\P_{u}(\tau^{d}<\infty)=(1-\mu)\sum_{j=u+1}^{\infty}p_{j}^{*(j-u)}-
\left(1-\P_1(\tau^{d}<\infty)\right)\sum\limits_{z=0}^{d-1}\P_{u}(\tau^0<\infty,-R_{\tau^0}=z)\P(\tau_{z+1}\leq d),\\
\end{eqnarray}
where
\begin{eqnarray}\label{foru=1}
\P_1(\tau^{d}<\infty)=\frac{\P_{1}(\tau^0<\infty)-\sum\limits_{z=0}^{d-1}\P_{1}(\tau^0<\infty,-R_{\tau^0}=z)\P(\tau_{z+1}\leq d)}{1-\sum\limits_{z=0}^{d-1}\P_{1}(\tau^0<\infty,-R_{\tau^0}=z)\P(\tau_{z+1}\leq d)}.
\end{eqnarray}
\end{Thm}
\proof
By taking limit as $t\rightarrow \infty$ in formula (\ref{finite}), we derive:
\begin{eqnarray}\nonumber
\P_{u}(\tau^{d}<\infty)=\P_u(\tau^0<\infty)
\left(1-\P_1(\tau^{d}<\infty)\right)\sum\limits_{s=1}^{\infty}\sum\limits_{\omega=1}^{d}\sum\limits_{z=0}^{\omega-1}\P_{u}(\tau^0=s,-R_{\tau^0}=z)\P(\tau_{z+1}=\omega),
\end{eqnarray}
where
$\P_u(\tau^0<\infty)$ is given in Lemma \ref{infty_Loisel}.
To get (\ref{infinite}) we use the following identity:
\begin{eqnarray}\nonumber
\sum\limits_{z=0}^{d-1}\P_{u}(\tau^0<\infty,-R_{\tau^0}=z)\P(\tau_{z+1}\leq d)=
\sum\limits_{s=1}^{\infty}\sum\limits_{\omega=1}^{d}\sum\limits_{z=0}^{\omega-1}\P_{u}(\tau^0=s,-R_{\tau^0}=z)\P(\tau_{z+1}=\omega).
\end{eqnarray}
To obtain (\ref{foru=1}) we apply (\ref{infinite}) with $u=1$.
\exit

\begin{Rem}\label{geometric2}\rm
Note that by Remark \ref{geometric}, in the compound binomial model with geometrical claim sizes,
the Parisian ruin probability can be found explicitly.
\end{Rem}

\section{Cram\'er's estimate of the ultimate Parisian ruin probability}\label{sec:cram}
In this section we derive the exponential asymptotics of the ultimate Parisian ruin probability
when a generic claim size has light-tailed distribution.

For all $\beta \geq 0$ we define moment-generating function
$$\varphi(\beta):=\log \E(e^{\beta R_1})=\beta +\log\left(\E(e^{-\beta Y})\right),$$
where $Y$ is a generic claim size and $\E_u$ is an expectation with respect to $\P_u$ (we skip subscript if $u=0$).
We will consider also a dual random walk $\widehat{R}_n=-R_n$ with a generic increment $U=Y-1$. Note that
$\P_u(\tau^0<\infty)=\P(\max_{k\geq 1} \widehat{R}_n>u)$.
Let $\widehat{L}^{-1}_n$ be the number of times new maxima are reached within $n$ steps of $\widehat{R}$. Let $\widehat{L}_n=\inf\{k\geq 1: \widehat{L}_k^{-1}=n \}$
be the number of steps required to achieve $n$ new maxima and $\widehat{H}_n=\widehat{R}_{\widehat{L}_n}$ be the $n$-th new maximum of $\widehat{R}$.
In other words, let $\{\widehat{L}_n, \widehat{H}_n), n\in \N\}$ be a ladder height process of $\widehat{R}$.

Assume Cram\'{e}r conditions that is: there exists $\gamma>0$ satisfying:
\begin{equation}\label{Cramerequation}
\varphi(-\gamma)=0
\end{equation}
and
\begin{equation}\label{Cramerequation2}
\varphi^\prime(-\gamma)<\infty.
\end{equation}

Above assumptions mean that $\E e^{\gamma Y}<\infty$ for a generic claim size $Y$, hence its distribution is light-tailed.

In lemmas below we recall the Cram\'er asymptotics for the ultimate classical ruin probability. Results follow from random walk theory
and renewal theory and as such seem to be classical ones. In the ruin theory for discrete risk process such results are presented e.g., in
Landriault et al. \cite{d} (see also Willmot and Lin \cite{Willmot2} and Rolski et al. \cite[p. 255-259]{Rolski}).
However, there exist different representation of the constant $C$ used in this results. In most cases,
these expressions for $C$ contain the ruin probability at $0$. Moreover authors often assume that claim size distribution is non-arithmetic,
which is not true in our case (since the span of the distribution equals one). For completeness we decided to present a proof based on Asmussen
\cite[Th. 13.5.2 and 13.5.3, p. 365]{asmappliedprob} (in a lattice version).

\begin{Fa}\label{Feller}
If we assume Cram\'{e}r conditions (\ref{Cramerequation}) and (\ref{Cramerequation2}), we have:
\begin{eqnarray} \label{cas}
\lim_{u \to \infty, u\in \N}e^{\gamma u} \P_u(\tau^0<\infty)= C
\end{eqnarray}
for
\begin{equation}\label{Ccram}
C=\frac{1-\mu }{-\varphi^\prime(-\gamma)}.\end{equation}
\end{Fa}
\proof
By Asmussen
\cite[Th. 13.5.3, p. 365]{asmappliedprob} it suffices to prove that constant (\ref{Ccram}) equals to:
\begin{equation}\label{stalaC}
C=\frac{1-\P(\widehat{L}_1<\infty)}{(1-e^{-\gamma})
\E \left[\widehat{H}_1 e^{\gamma \widehat{H}_1}, \widehat{L}_1<\infty)\right]}.
\end{equation}
For process $R$ we can define weakly ascending height process $\{(L_n, H_n), n\in \N\}$.
Note that $H_1=(1-p_0)\delta_0+p_0\delta_1$, where $\delta_x$ denotes the Dirac's Delta at $x$.
From the Wiener-Hopf factorization (see Asmussen \cite[ p. 234]{asmappliedprob}):
\begin{equation}\label{baxter2}
1-\E \left[e^{\theta \widehat{H}_1}, \widehat{L}_1<\infty\right]=\frac{1-e^{\varphi(-\theta)}}{(1-e^{-\theta})p_0}.
\end{equation}
Moreover, from (\ref{baxter2}) when $\theta \searrow 0$ we have
\begin{equation}\label{baxter}
1-\P(\widehat{L}_1<\infty)=\varphi'(0+)=(1-\mu)/p_0.
\end{equation}
To complete the proof we take derivative at $\theta=\gamma$ of the equation (\ref{baxter2}). Then we plug its result and (\ref{baxter}) into the right-hand side of the equation (\ref{stalaC}).
\exit

\begin{Fa}\label{Feller_parisian}
If we assume Cram\'{e}r condition (\ref{Cramerequation}), for $\theta>\gamma$
we have:
\begin{eqnarray} \label{cas_parisian}
\lim_{u \to \infty, u\in \N}\E_u\left[e^{\theta R_{\tau^0}}|\tau^0<\infty\right]= D(\theta),
\end{eqnarray}
where $$D(\theta)=\frac{(1-e^{\varphi(\theta)})(1-e^{-\gamma})}{(1-\mu)(1-e^{\theta})(1-e^{-(\gamma+\theta)})}.$$
\end{Fa}
\proof
Observe first that:
$$\E_u\left[e^{\theta R_{\tau^0}}|\tau^0<\infty\right]=\E\left[ e^{-\theta (\widehat{R}_{\widehat{\tau}_u}-u)}|\widehat{\tau}_u<\infty\right],
$$
where $\widehat{\tau}_u=\inf\{k\in \N: \widehat{R}_k \geq u\}$.
Define Esscher transform via:
\begin{equation}\label{Girsanov}
\left. \frac{d\P_{u}^{c}}{d\P_{u}}\right| _{\mathcal{F}_{n}}=\frac{\mathcal{E}%
_{n}\left( c\right) }{\mathcal{E}_{0}\left( c\right) }
\end{equation}
for any $c$ such that $\E e^{c \widehat{R}_1}<\infty$, where $\mathcal{E}_{n}\left( c\right) =\exp \{c\widehat{R}_{n}-\varphi
\left( -c\right) n\}$ is the exponential martingale under $\P_u$.
From Asmussen \cite[Th.13.5.2 and 13.5.3, p. 365]{asmappliedprob}
we have:
$$\lim_{u\to\infty, u\in \N}\E\left[e^{-\theta (\widehat{R}_{\widehat{\tau}_u}-u)}|\widehat{\tau}_u<\infty\right]=
\lim_{u\to\infty, u\in \N}\E^{\gamma}\left[e^{-(\theta+\gamma) (\widehat{R}_{\widehat{\tau}_u}-u)}\right]/\E^{\gamma}\left[e^{-\gamma (\widehat{R}_{\widehat{\tau}_u}-u)}\right].$$
Finally, from Asmussen \cite[Th. VIII.2.1, p. 224]{asmappliedprob} we obtain:
\begin{eqnarray}\nonumber
\lim_{u\to\infty, u\in \N}\E^{\gamma}\left[e^{-\theta (\widehat{R}_{\widehat{\tau}_u}-u)}\right]=
\frac{1}{\E^\gamma \widehat{H}_1}\sum_{k=0}^\infty e^{-\theta k}(1-\P^\gamma(\widehat{H}_1\leq k))
=\frac{1-\E \left[e^{(\gamma-\theta)\widehat{H}_1}, \widehat{L}_1<\infty\right]}{(1-e^{-\theta})\E \left[ e^{\gamma \widehat{H}_1} \widehat{H}_1, \widehat{L}_1<\infty\right]}.\label{falka1}\\
\end{eqnarray}\\
Applying (\ref{baxter2}) to (\ref{falka1}) completes the proof.
\exit\\

Let $\Phi(\cdot)$ be the inverse of $\varphi$.
We are ready now to state the main result of this section.
\begin{Thm}\label{Crameras}
Under the Cram\'er conditions (\ref{Cramerequation})- (\ref{Cramerequation2}) we have:
\begin{eqnarray} \label{cas_par}
\lim_{u \to \infty, u\in \N}e^{\gamma u} \P_u(\tau^d<\infty)
&=& C\left[1-(1-\P_1 (\tau^d<\infty))f(d)\right],
\end{eqnarray}
with $C$ defined in (\ref{stalaC}), $\P_1(\tau^d<\infty)$ given in Theorem \ref{ThmMain2},
where, for $\theta>\varphi(\gamma)$:
\begin{equation}\label{LTfk}
\sum_{k=0}^{\infty}e^{-k\theta} f(k)=
\frac{e^{-\Phi(\theta)}}{(1-e^{-\theta})}D(\Phi(\theta)).
\end{equation}
\end{Thm}

\proof
To prove Cram\'er asymptotics (\ref{cas_par}) we use Lemma \ref{infty_Loisel}, equations (\ref{infinite}) and (\ref{cas}), where
\begin{eqnarray}\label{sumy}\nonumber
f(k)=\sum\limits_{z=0}^{k-1}\lim_{u\to\infty, u\in \N}\P_{u}(-R_{\tau^0}=z|\tau^0<\infty)\P(\tau_{z+1}\leq d).
\end{eqnarray}
To prove that function $f$ has the Laplace transform (\ref{LTfk}) observe that:
 \begin{eqnarray}\nonumber
\sum_{k=0}^{\infty}e^{-k\theta}\P(\tau_{z+1}\leq k)=\frac{1}{(1-e^{-\theta})}\sum_{k=0}^{\infty}e^{-(k+1)\theta}\P(\tau_{z+1}= k+1)=\frac{\E(e^{-\theta\tau_{z+1}},\tau_{z+1}<\infty)}{(1-e^{-\theta})}=\frac{e^{-\Phi(\theta)(z+1)}}{(1-e^{-\theta})}.
\end{eqnarray}
In the first equality we use summation-by-parts formula. The last equality is a consequence of
Optional Stopping Theorem applied to the martingale $\mathcal{E}_{n}\left( \Phi(\theta)\right)$ at the stopping time $\tau_z$.
Hence, by Lemma \ref{Feller_parisian}:
\begin{eqnarray}\nonumber
\sum_{k=0}^{\infty}e^{-k\theta}f(k)
= \frac{e^{-\Phi(\theta)}}{(1-e^{-\theta})}\lim_{u\to\infty, u\in \N}\E_{u}\left(e^{\Phi(\theta)R_{\tau^0}}|\tau^0<\infty\right)=
 \frac{e^{-\Phi(\theta)}}{(1-e^{-\theta})}D(\Phi(\theta)).
\end{eqnarray}
\exit

\section{Heavy-tailed estimate of the ultimate Parisian ruin probability}\label{sec:heavy}

In this section we will assume that Cram\'er equation
(\ref{Cramerequation}) has no solution.
In particular, we assume that the distribution $\{p_n, n\in \N\}$ of $Y$
belongs to the class $\mathcal{S}^{(\alpha)}$.
We refer to Asmussen and Albrecher \cite{Assbook} and Foss et al. \cite{Serdimastan} for all properties of these class of distributions;
see also Tang et al. \cite{Tang1, Tang2} and references therein.
This class is defined as follows.

\begin{Def} (Class $\mathcal{L}^{(\alpha)}$, lattice case of span $1$)
For a parameter $\alpha \geq 0$ we say that distribution function $G$ on $\N$ with tail $\overline{G}=1-G$
belongs to class $\mathcal{L}^{(\alpha)}$ if
\begin{itemize}
\item[$(i)$]  $\overline{G}(k)>0$ for each $k \in \N$,
\item[$(ii)$] $\lim_{n \rightarrow \infty} \frac{\overline{G}(n-1)}{\overline{G}(n)}=e^{\alpha}.$
\end{itemize}

\end{Def}

\begin{Def}\label{def2} (Class $\mathcal{S}^{(\alpha)}$, lattice case of span $1$)
We say that $G$ belongs to class $\mathcal{S}^{(\alpha)}$
if
\begin{itemize}
\item[$(i)$] $G \in \mathcal{L}^{(\alpha)}$,
\item[$(ii)$] $\sum_{k=0}^\infty e^{\alpha k} (G(k)-G(k-1))<\infty$,
\item[$(iii)$] for some $M_0<\infty$, we have
\begin{eqnarray}
\lim_{u \rightarrow \infty, u\in\N} \frac{\overline{G^{*2}}(u)}{\overline{G}(u)}=2M_0, \label{convvv}
\end{eqnarray}
where $\overline{G^{*2}}(u)=1-G^{*2}(u)$ and $*$ denotes convolution.
\end{itemize}
\end{Def}

If $G\in \mathcal{S}^{(\alpha)}$ then we say that $G$ is convolution equivalent.
The case $\alpha=0$ is particularly interesting since the class $\mathcal{S}^{(0)}$ is a class
of subexponential distributions. If $G\in \mathcal{S}^{(0)}$ then
it is heavy-tailed and its moment generating function does not exist for any
strictly positive arguments and therefore right-hand side of (\ref{Cramerequation}) is not well-defined.
Distributions with regularly varying tails are also in class  $\mathcal{S}^{(0)}$. Typical example is the Pareto distribution.


Let $F(k)=\sum_{l=0}^kp_l$ and $F_I(k)=\frac{1}{\mu}\sum_{l=0}^k\overline{F}(l)$, $k\in \N$.
From now on we will assume that
either
\begin{equation}\label{assub1}
F_I\in \mathcal{S}^{(0)}
\end{equation}
or
\begin{equation}\label{assub2}
F\in \mathcal{S}^{(\alpha)},\qquad \alpha>0.
\end{equation}

We recall now the asymptotic result for the ultimate ruin probability.
\begin{Fa}\label{sub}
If (\ref{assub1}) holds then
$$\lim_{u\to\infty, u\in \N}\frac{\P_u(\tau^0<\infty)}{\frac{\mu}{1-\mu} \overline{F}_I(u)}=1.$$
If (\ref{assub2}) holds then
$$\lim_{u\to\infty, u\in \N}\frac{\P_u(\tau^0<\infty)}{K\overline{F}(u)}=1,$$
where
\begin{equation}\label{K}
K=\frac{(1-\mu)(1-e^{-\alpha})}{(1-e^{\varphi(-\alpha)})^2}.\end{equation}
\end{Fa}
\proof
For the case $\alpha =0$ observe that the ruin probability $\P_u(\tau^0<\infty)$ equals the
ruin probability for the classical renewal risk process in a continuous time with interarrival time
equal $1$ and the generic claim size $Y$. Then the first part of assertion follows from Asmussen and Albrecher \cite[Th. 10.3.1, p. 305]{Assbook} (see also Foss et al. \cite[Th. 5.12, p. 113]{Serdimastan}).
The case $\alpha >0$ follows from Bertoin and Doney \cite[Th. 1]{BertDon2} and Asmussen \cite[eq. (4.4.5), p. 231]{asmappliedprob} with
$$K=\frac{1-\P(\widehat{L}_1<\infty)}{(1-e^{\varphi(-\alpha)})(1-\E[ e^{\alpha \widehat{H}_1}, \widehat{L}_1<\infty])}.$$
Identities (\ref{baxter}) and (\ref{baxter2}) complete the proof.
\exit

Observe that
$$\P_u(\tau^0<\infty)=\P(\tau(u)<\infty)$$
and that
$$\P_u(R_{\tau^0}=k|\tau^0<\infty)=\P(X_{\tau(u)}=k|\tau(u)<\infty),\qquad k\in \N,$$
where
$$\tau(u)=\inf\{t\geq 0: X_t>u\}$$
for a compound Poisson process $X_t=\sum_{i=1}^{N_t}(Y_i-1)$ and
$N_t$ being independent of $\{Y_i\}$ Poisson process with intensity $1$.
It is enough to observe this continuous time process $X_t$ at the moments
of jumps. We transfer discrete risk process into continuous-time set-up just
to use almost straightforward convenient references concerning L\'{e}vy processes
and to avoid more direct and longer proofs.
From Kl$\ddot{\textrm{u}}$ppelberg and Kyprianou \cite[Th. 4.2, eq. (2.8) and Remark 4.3 (iii)]{KyprKlup2} we have that
\begin{Fa}\label{andreas}
There exists function $\overline{W}$ such that:
\begin{equation*}
\lim_{u \to \infty, u\in \N}\P(X_{\tau(u)}\geq k|\tau(u)<\infty)=\overline{W}(k).
\end{equation*}
If (\ref{assub1}) holds then the distribution $W=1-\overline{W}$ is degenerate placing
all its mass at infinity which follows from so-called principle of one big jump.
If (\ref{assub2}) holds then the function
$\overline{W}$ is a tail of (possibly improper) distribution function:
\begin{equation}\label{W}
\overline{W}(k)=\frac{e^{-\alpha k}}{1-\mu}\left(\frac{-\phi(-\alpha)}{\alpha}+\sum_{l=k+1}^\infty
\left(e^{\alpha l}-e^{\alpha k}\right)\P(\widehat{H}_1=l, \widehat{L}_1<\infty)\right)
\end{equation}
with $\phi(\alpha)=\log\E e^{\alpha X_1}=e^{-\alpha}\E e^{\alpha Y}- 1$.
\end{Fa}

\begin{Rem}\rm
\rm Note that by Definition \ref{def2}(ii) for $\alpha >0$
we have $\E e^{\alpha Y}<\infty$ and $\E \left[e^{\alpha \widehat{H}_1}, \widehat{L}_1<\infty\right]<\infty$.
Unfortunately, it seems that function $\overline{W}$ is very hard to identify more explicitly.
\end{Rem}

The main result of this section is following.
\begin{Thm}\label{subpar}
Under assumption (\ref{assub1}) we have:
$$\lim_{u \to \infty, u\in \N}\frac{\P_u(\tau^d<\infty)}{\overline{F}_I(u)}
= \frac{\mu}{1-\mu} .
$$
Under assumption (\ref{assub2}) we have:
$$\lim_{u \to \infty, u\in \N}\frac{\P_u(\tau^d<\infty)}{\overline{F}(u)}=BK,$$
where constant $K$ is given in (\ref{K}) and $B=
\left[1-(1-\P_1 (\tau^d<\infty))g(d)\right]$ for a function $g$ on $\N$ with the Laplace transform:
\begin{equation}\label{LTfk2}
\sum_{k=0}^{\infty}e^{-k\theta} g(k)=
\frac{e^{-\Phi(\theta)}}{(1-e^{-\theta})}w(\Phi(\theta)),
\end{equation}
where
$$w(\theta)=\frac{e^{-\Phi(\theta)}}{(1-e^{-\theta})}
\sum_{k=1}^\infty e^{-\Phi(\theta)k}\left(\overline{W}(k-1)-\overline{W}(k)\right).$$
\end{Thm}
\proof
The proof is similar to the proof of Theorem \ref{Crameras}, where we use
Lemma \ref{andreas} instead of Lemma \ref{Feller}.
\exit

\begin{Rem}\rm
\rm
Note that the Parisian delay has influence on the heavy-tailed asymptotics of the ruin probability only when $\alpha >0$. The subexponential case (when $\alpha=0$) gives the same asymptotics as for classical ruin moment $\tau^0$, which is a consequence of one big claim that causes the ruin.
\end{Rem}

\newpage
\section{Examples}\label{sec:examples}

In this section we use Theorem \ref{ThmMain} to calculate the Parisian non-ruin probability
for various initial capitals and Parisian delays.
We consider classical binomial model (\ref{CL}), where generic claim is a product of Bernoulli random variable $I$ with
$\P(I=1)=1-\P(I=0)=b>0$ and
some other positive random variable $B$ with values in natural numbers.
To capture different behaviors of the Parisian non-ruin probability, we consider light-tailed case (in Example 1 with geometric d.f. of the claim size)
and heavy-tailed case (in Example 2 with Pareto d.f. of the claim size). In the first case we also found the ultimate non-ruin probability $P_{u}(\tau^d=\infty)$
and compared it with the finite-time one. All numerical calculations have been made in the R package.\\

{\bf Example 1.}

In this example we consider generic claim $Y=I\cdot B$ being the product of Bernoulli random variable $I$ with
$\P(I=1)=1-\P(I=0)=b>0$ and geometric random variable $B$ with parameter $1-q>0$
giving the distribution of the claim size (\ref{bingeom}).
We assume that $b=0.08$ and $1-q=0.1$. Hence the claim size mean equals $\mu=\E Y_1=0.8<1$ which gives positive safety loading.

\begin{figure}[h!]
\caption{Parisian finite-time and ultimate non-ruin probability for different initial capitals - Geometric claims}
\label{rysGeometric2}
 \centering
  \includegraphics[width=1.1\textwidth]{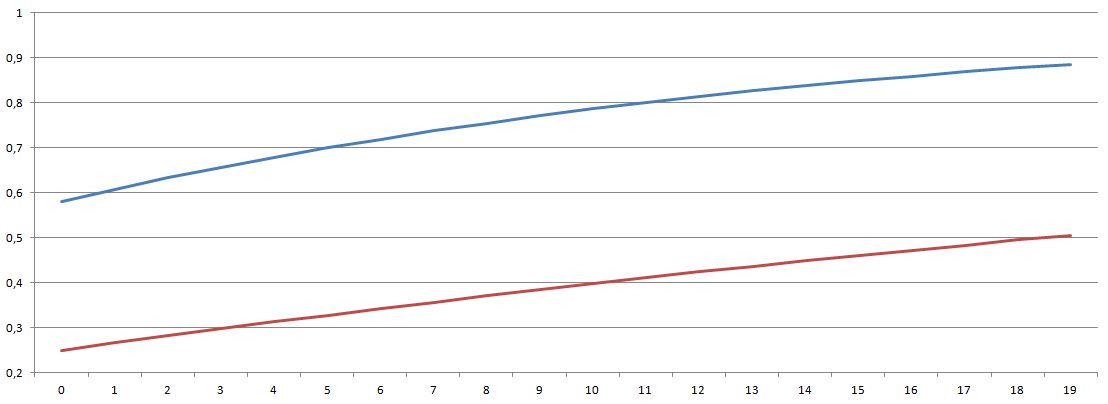}
\end{figure}

\begin{figure}[h!]
\caption{Parisian finite-time and ultimate non-ruin probability for different Parisian delay - Geometric claims}
\label{rysGeometric3}
 \centering
  \includegraphics[width=1.1\textwidth]{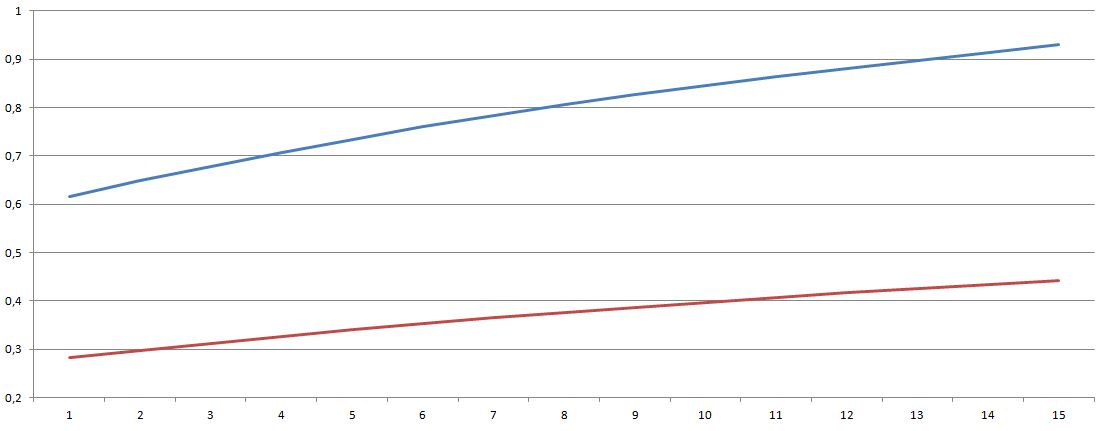}
\end{figure}

At the beginning we take $u=4$ (initial capital) and $d=3$ (Parisian delay). Table \ref{tabelaGeometric1} identifies Parisian non-ruin probability for different $t\leq 27$.
All calculations were performed using Theorem \ref{ThmMain}.
Similar result could be derived for the Parisian non-ruin probability for fixed time horizon $t=20$
with different initial capitals (see Table \ref{tabelaGeometric2} and Figure \ref{rysGeometric2}) and for different Parisian delays (see Table \ref{tabelaGeometric3} and Figure \ref{rysGeometric3}). For comparison we also added ultimate non-ruin probability $P_{u}(\tau^{d}=\infty)$ which was found using Remark \ref{geometric} and Theorem \ref{ThmMain2}. Note that the difference $\Delta$ between these two quantities gives the probability of ruin after or at time $t=20$. In Figures \ref{rysGeometric2} and \ref{rysGeometric3} we use red color to denote ultimate non-ruin probability, and blue color to denote the non-ruin probability over finite-time horizon.

All these calculations show that the formula given in Theorem \ref{ThmMain}
produces deep comparison results for very wide choice of parameters.
In particular, they demonstrate that increasing the Parisian delay can substantially decrease the ruin probability, while keeping initial capital fixed.
Moreover, considering the difference $\Delta=\P_u(\tau^d\geq 20)-\P_u(\tau^d =\infty)$, it seems plausible that the ruin will happen after time $t=20$ and hence after long time evolution of the risk process (\ref{CL}).\\

{\bf Example 2.}

To analyze the heavy-tailed case
in this the example, we consider
generic claim $Y=I\cdot B$ being the product of Bernoulli random variable $I$ with
$\P(I=1)=1-\P(I=0)=b=0.08$ and Pareto random variable $B$. That is, we have:
\begin{equation}\label{Pareto}\nonumber
\P(Y=0)=1-b,\qquad \P(Y=z)=b\left(\frac{1}{z^\alpha}-\frac{1}{(z+1)^\alpha}\right) \quad \text{for $z=1,2,\ldots .$}
\end{equation}
Note that then $\mathbb{E}Y=b\zeta(\alpha)$, where $\zeta$ is Riemann zeta function. To obtain $\mathbb{E}Y=0.8$ as it was in the previous example, we take $\alpha = 1.1062123$.

\begin{figure}[h!]
\caption{Parisian non-ruin probability for different time horizons - Pareto claims}
\label{rysPareto1}
 \centering
  \includegraphics[width=1.1\textwidth]{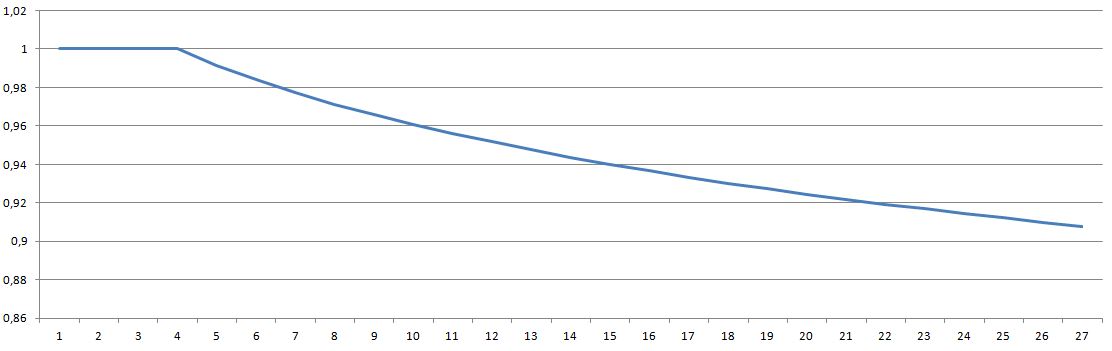}
\end{figure}

\begin{figure}[h!]
\caption{Parisian non-ruin probability for different initial capitals and $t=20$ - Pareto claims}
\label{rysPareto2}
 \centering
  \includegraphics[width=1.1\textwidth]{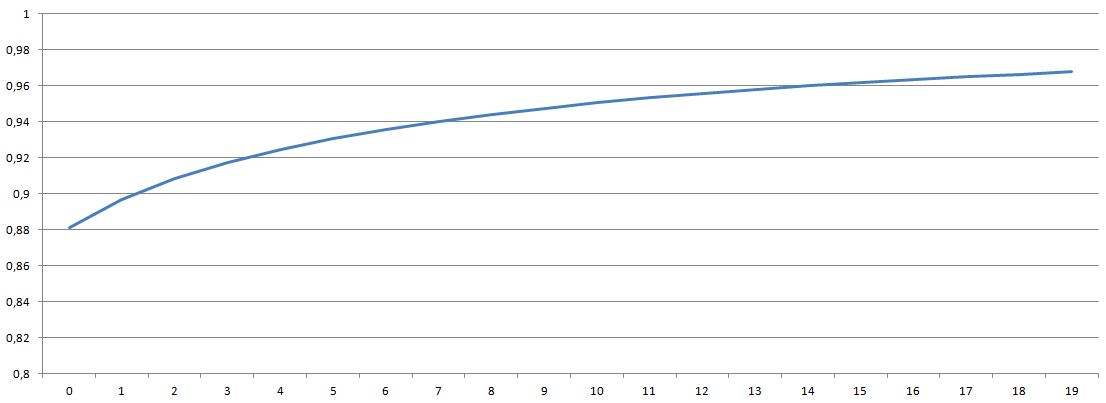}
\end{figure}

\begin{figure}[h!]
\caption{Parisian non-ruin probability for different Parisian delay and $t=20$ - Pareto claims}
\label{rysPareto3}
 \centering
  \includegraphics[width=1.1\textwidth]{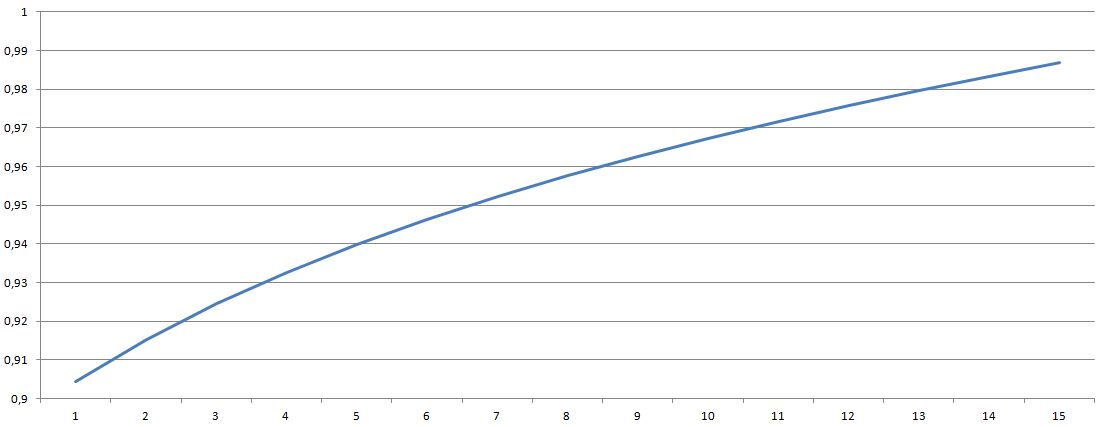}
\end{figure}

In the tables and figures below, we show how the Parisian non-ruin probability changes for different time horizons (Table \ref{tabelaPareto1} and Figure \ref{rysPareto1}
with $u=5$ (initial capital), $d=3$ (Parisian delay)), different initial capitals (Table \ref{tabelaPareto2} and Figure \ref{rysPareto2} with $d=3$ (Parisian delay) and $t=20$ (time horizon)) and different Parisian delays (Table \ref{tabelaPareto3} and Figure \ref{rysPareto3} with $u=5$ (initial capital) and $t=20$ (time horizon)).

In the heavy-tailed case, the non-ruin Parisian probability is much bigger than in the light-tailed case. At the same time in the heavy-tailed case, the loss given default is substantially larger.

\newpage

\begin{table}[t]
\caption{Parisian non-ruin probability $\P_{u=4}(\tau^{d=3}\geq t)$ for different time horizons - Geometric claims}
\label{tabelaGeometric1}
\begin{center}
\begin{tabular}{|*{2}{l|}} \hline
Time t & \textbf{$\P_{u=4}(\tau^{d=3}\geq t)$}\\ \hline
1..4 & 1 \\ \hline
5 & 0,959785 \\ \hline
6 & 0,925200 \\ \hline
7 & 0,894939 \\ \hline
8 & 0,868044 \\ \hline
9 & 0,843803 \\ \hline
10 & 0,821846 \\ \hline
11 & 0,801862 \\ \hline
12 & 0,783589 \\ \hline
13 & 0,766809 \\ \hline
14 & 0,751338 \\ \hline
15 & 0,737022 \\ \hline
16 & 0,723729 \\ \hline
17 & 0,711349 \\ \hline
18 & 0,699784 \\ \hline
19 & 0,688951 \\ \hline
\textbf{20} & \textbf{0,678780} \\ \hline
21 & 0,669207 \\ \hline
22 & 0,660177 \\ \hline
23 & 0,651642 \\ \hline
24 & 0,643560 \\ \hline
25 & 0,635894 \\ \hline
26 & 0,628609 \\ \hline
27 & 0,621676 \\ \hline
\end{tabular}
\end{center}
\end{table}

\begin{figure}[h!]
\caption{Parisian non-ruin probability $\P_{u=4}(\tau^{d=3}\geq t)$ for different time horizons - Geometric claims}
\label{rysGeometric1}
 \centering
  \includegraphics[width=1.1\textwidth]{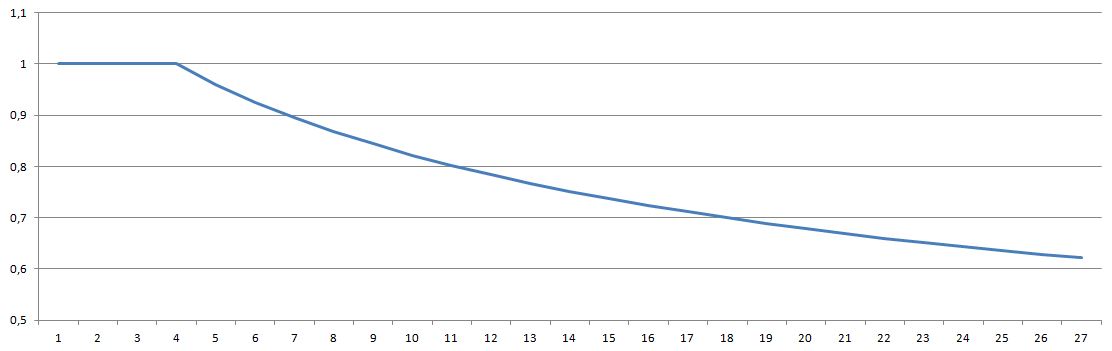}
\end{figure}

Results of table above were also presented on Figure \ref{rysGeometric1}.

\newpage
\begin{table}[t]
\caption{Parisian finite-time and ultimate non-ruin probability for different initial capitals - Geometric claims}
\label{tabelaGeometric2}
\begin{center}
\begin{tabular}{|*{3}{c|c|c|}} \hline
u & \textbf{ $\P_{u}(\tau^{d=3} \geq 20)$}& \textbf{$\P_{u}(\tau^{d=3}=\infty)$}&$\Delta=\P_{u}(\tau^{d=3} \geq 20)-\P_{u}(\tau^{d=3} =\infty )$ \\ \hline
0 & 0,5810479 & 0,249772		& 	0,3312759
																
\\ \hline
1 & 0,607774 & 0,266081 & 0,3416928\\ \hline
2 & 0,632917 & 0,282036 & 0,3508806\\ \hline
3 & 0,656559 & 0,297644 & 0,3589148\\ \hline
4 & \textbf{0,678780} & 0,312913 & 0,3658669\\ \hline
5 & 0,699656  &0,327849 & 0,3718072\\ \hline
6 & 0,719260  &0,342461 & 0,3767995\\ \hline
7 & 0,737663 & 0,356756 & 0,3809066\\ \hline
8 & 0,754929 & 0,370739 & 0,3841904\\ \hline
9 & 0,771124 & 0,384418 & 0,3867056\\ \hline
10 & 0,786308 &0,397801 & 0,3885074\\ \hline
11 & 0,800539 &0,410892 & 0,3896471\\ \hline
12 & 0,813871 &0,423699 & 0,3901725\\ \hline
13 & 0,826358 &0,436227 & 0,3901309\\ \hline
14 & 0,838048 &0,448483 & 0,3895649\\ \hline
15 & 0,848989 &0,460473 & 0,3885157\\ \hline
16 & 0,859225 &0,472202 & 0,3870228\\ \hline
17 & 0,868799 &0,494899 & 0,3738986\\ \hline
18 & 0,877750 &0,483675 & 0,3940752\\ \hline
19 & 0,886117 &0,505880 & 0,3802373 \\ \hline
\end{tabular}
\end{center}
\end{table}

\newpage

\begin{table}[t]
\caption{Parisian finite-time and ultimate non-ruin probability for different Parisian delays - Geometric claims}
\label{tabelaGeometric3}
\begin{center}
\begin{tabular}{|*{3}{c|c|c|}} \hline
\textbf{$d$} & $P_{u=4}(\tau^{d} \geq 20)$& $P_{u=4}(\tau^{d} =\infty )$&$\Delta=P_{u=4}(\tau^{d} \geq 20)-P_{u=4}(\tau^{d} =\infty )$ \\ \hline
1 & 0,615985 & 0,283120 & 0,332865\\ \hline
2 & 0,648228 & 0,298331 & 0,349897	\\ \hline
3 & \textbf{0,678780} & 0,312913& 0,365867\\ \hline
4 & 0,707581 & 0,326841 &	0,380740\\ \hline
5 & 0,734634 & 0,340117 & 0,394517	\\ \hline
6 & 0,759986 & 0,352754 &	0,407232\\ \hline
7 & 0,783716 & 0,364778 & 0,418938	\\ \hline
8 & 0,805913 & 0,376220	& 0,429693\\ \hline
9 & 0,826625 & 0,387117	& 0,439508\\ \hline
10 & 0,845859 & 0,397502 & 0,448357	\\ \hline
11 & 0,863890 & 0,407412	& 0,456479\\ \hline
12 & 0,881019 & 0,416880 &	0,464139\\ \hline
13 & 0,897518 & 0,425939 & 0,471579\\ \hline
14 & 0,913656 & 0,434617 & 0,479039	\\ \hline
15 & 0,929708 & 0,442944  & 0,486764\\ \hline
\end{tabular}
\end{center}
\end{table}

\newpage

\begin{table}[t]
\caption{Parisian non-ruin probability for different time horizons - Pareto claims}
\label{tabelaPareto1}
\begin{center}
\begin{tabular}{|*{2}{l|}} \hline
Time t & \textbf{$\P_{u=4}(\tau^{d=3}\geq t)$}\\ \hline
1..4 & 1 \\ \hline
5 & 0,991491 \\ \hline
6 & 0,984043 \\ \hline
7 & 0,977390 \\ \hline
8 & 0,971360 \\ \hline
9 & 0,965837 \\ \hline
10 & 0,960746 \\ \hline
11 & 0,956030 \\ \hline
12 & 0,951638 \\ \hline
13 & 0,947532 \\ \hline
14 & 0,943676 \\ \hline
15 & 0,940047 \\ \hline
16 & 0,936617 \\ \hline
17 & 0,933368 \\ \hline
18 & 0,930281 \\ \hline
19 & 0,927343 \\ \hline
20 & \textbf{0,924540} \\ \hline
21 & 0,921860 \\ \hline
22 & 0,919294 \\ \hline
23 & 0,916834 \\ \hline
24 & 0,914470 \\ \hline
25 & 0,912195 \\ \hline
26 & 0,910005 \\ \hline
27 & 0,907892 \\ \hline
\end{tabular}
\end{center}
\end{table}

\newpage

\begin{table}[t]
\caption{Parisian non-ruin probability for different initial capitals and $t=20$ - Pareto claims}
\label{tabelaPareto2}
\begin{center}
\begin{tabular}{|*{2}{l|}} \hline
u & \textbf{$P_{u}(\tau^{d=3} \geq 20)$}\\ \hline
0 & 0,881454 \\ \hline
1 & 0,896836 \\ \hline
2 & 0,908254 \\ \hline
3 & 0,917233 \\ \hline
4 & \textbf{0,924540} \\ \hline
5 & 0,930631 \\ \hline
6 & 0,935802 \\ \hline
7 & 0,940255 \\ \hline
8 & 0,944135 \\ \hline
9 & 0,947548 \\ \hline
10 & 0,950576 \\ \hline
11 & 0,953289 \\ \hline
12 & 0,955714 \\ \hline
13 & 0,957914 \\ \hline
14 & 0,959912 \\ \hline
15 & 0,961735 \\ \hline
16 & 0,963406 \\ \hline
17 & 0,964943 \\ \hline
18 & 0,966360 \\ \hline
19 & 0,967673 \\ \hline
\end{tabular}
\end{center}
\end{table}

\newpage
\begin{table}[t]
\caption{Parisian non-ruin probability for different Parisian delays and $t=20$ - Pareto claims}
\label{tabelaPareto3}
\begin{center}
\begin{tabular}{|*{2}{l|}} \hline
\textbf{$d$} & \textbf{$P_{u=4}(\tau^{d} \geq 20)$}\\ \hline
1 & 0,904499 \\ \hline
2 & 0,915302 \\ \hline
3 & \textbf{0,92454} \\ \hline
4 & 0,932625 \\ \hline
5 & 0,939821 \\ \hline
6 & 0,946308 \\ \hline
7 & 0,952214 \\ \hline
8 & 0,957633 \\ \hline
9 & 0,962638 \\ \hline
10 & 0,967283 \\ \hline
11 & 0,971624 \\ \hline
12 & 0,975709 \\ \hline
13 & 0,979579 \\ \hline
14 & 0,983266 \\ \hline
15 & 0,986801 \\ \hline
\end{tabular}
\end{center}
\end{table}

\end{document}